%% file: main.tex
\documentclass[conference,10pt,letterpaper]{IEEEtran}

\usepackage{amsmath,amssymb,amsfonts,amsthm}
\usepackage{graphicx}
\usepackage{xcolor}
\usepackage{bm}
\usepackage[caption=false,font=footnotesize]{subfig}
\usepackage{algorithm}
\usepackage{algorithmicx}
\usepackage[noend]{algpseudocode}
\definecolor{commBdy}{HTML}{D55E00} % boundary communication
\definecolor{commLvl}{HTML}{0072B2} % level-transition communication
\newcommand{\algcmt}[1]{\quad\(\triangleright\)~\textit{#1}}

\usepackage{url}
\usepackage{booktabs}
\usepackage{color, colortbl}
\usepackage{mathtools}
\usepackage{multirow}
\usepackage{siunitx}
\usepackage{epstopdf}
\usepackage{cite}
\usepackage[hidelinks]{hyperref}
\usepackage{cleveref}
\usepackage{tikz}
\usetikzlibrary{arrows.meta,calc}

\theoremstyle{definition}

\theoremstyle{remark}

\newcommand{\bigo}[1]{\mathcal{O}\left(#1 \right) }

\input{notations}

\ifdefined\IPDPSsubmission
  \newcommand{\chao}[1]{}
\else
  \newcommand{\chao}[1]{{\color{red}[ Chao: #1 ]}}
\fi

\title{Scaling $\mathcal{H}^2$ Factorization to Billion-Unknown Three-Dimensional Problems}

\author{
\IEEEauthorblockN{Tianyu Liang\IEEEauthorrefmark{1},
Chao Chen\IEEEauthorrefmark{2},
James Demmel\IEEEauthorrefmark{1},
and Yang Liu\IEEEauthorrefmark{3}}
\IEEEauthorblockA{\IEEEauthorrefmark{1}UC Berkeley\qquad
\IEEEauthorrefmark{2}NC State University\qquad
\IEEEauthorrefmark{3}Lawrence Berkeley National Laboratory}
\IEEEauthorblockA{\texttt{tianyul@berkeley.edu}, \texttt{chao\_chen@ncsu.edu},\\
\texttt{demmel@berkeley.edu}, \texttt{liuyangzhuan@lbl.gov}}
}

\begin{document}
\bstctlcite{IPDPS2027:BSTcontrol}

\maketitle

\begin{abstract}
\input{abstract}
\end{abstract}


\input{introduction}

\input{background}

\input{distributed-algorithm}

\input{lazy-schur}

\input{streaming-id}

\input{shared-memory}

\input{analysis}

\input{experiments}

% Old numerical-results section from the previous paper.  Keep results.tex as a
% reference for table/figure style, but exclude it from this draft.
%\input{results}

\input{conclusion}

\bibliographystyle{IEEEtran}
\bibliography{biblio}

\end{document}

%% file: notations.tex
\newcommand{\N}{ {\mathcal{N}} }

\newcommand{\bigO}{ {\mathcal{O}} }
\newcommand{\B}{ {\mathcal{B}} }
\newcommand{\X}{ {\mathcal{X}} }
\newcommand{\F}{ {\mathcal{F}} }
\newcommand{\T}{ {\mathcal{T}} }
\newcommand{\eps}{ {\varepsilon} }
\newcommand{\I}{ {\mathcal{I}} }
\newcommand{\J}{ {\mathcal{J}} }
\newcommand{\R}{ {\mathcal{R}} }
\newcommand{\C}{ {\mathcal{C}} }
\newcommand{\M}{ {\mathcal{M}} }

\renewcommand{\P}{ {\mathcal{P}} }
\renewcommand{\S}{ {\mathcal{S}} }

\newcolumntype{H}{>{\setbox0=\hbox\bgroup}c<{\egroup}@{}}

\definecolor{Gray}{gray}{0.9}

%% file: abstract.tex
Strong recursive skeletonization provides approximate direct solvers for dense matrices from integral equations and kernel methods, but distributed implementations face substantial communication costs and limited parallelism, especially in three dimensions. We present a distributed-memory $\mathcal{H}^2$ factorization solver for symmetric matrices that combines spatial box coloring, delayed Schur complement updates, streamed sketching-based interpolative decomposition, and BLAS-3 and shared-memory kernels. Spatial coloring enables every process to participate in each boundary stage, while delayed updates and streamed sketching reduce dense storage and communication. To our knowledge, this is the first distributed-memory strong recursive skeletonization factorization demonstrated for three-dimensional volume integral equations at billion-unknown scale.

Experiments on up to 512 compute nodes demonstrate factorization of a three-dimensional Laplace system with 3.62 billion unknowns, alongside three-dimensional Helmholtz tests and a two-dimensional Mat\'ern system with 268 million unknowns. Combining streamed sketching and fully delayed updates reduces factorization time by $6.9\times$ for Laplace and $8.5\times$ for Helmholtz relative to the unsketched, eager configuration, while reducing maximum per-process resident memory by 57\% and 54\%, respectively.

%% file: introduction.tex
\section{Introduction}

Dense linear systems arise in integral-equation models of electrostatics,
wave propagation, and fluid flow, and in covariance-based statistical
models. Conventional dense factorization requires $\bigo{N^3}$ arithmetic
and $\bigo{N^2}$ storage for $N$ unknowns. These costs preclude the use of direct solvers
at the scales required by modern simulations. Rank-structured solvers
exploit low-rank off-diagonal interactions to reduce these costs to quasi-linear or linear complexity \cite{martinsson2025fast,li2026parallel}. Such a factorization can serve as an approximate direct solver for multiple right-hand sides or as a preconditioner for iterative solution when ill-conditioning makes an unpreconditioned method ineffective.

We study strongly admissible $\mathcal{H}^2$ factorizations of symmetric matrices $A=A^\top$, including real-symmetric and complex-symmetric kernels. Throughout, $\top$ denotes transpose without conjugation; positive definiteness is not assumed. These factorizations compress
interactions between well-separated clusters and use nested bases across
a spatial hierarchy \cite{hackbusch2003introduction,hackbusch2002data,minden2017recursive}.
Strong recursive skeletonization (RS-S) \cite{minden2017recursive}, among several other $\mathcal{H}^2$ factorization variants such as FMM-LU
\cite{sushnikova2023fmm}, inverse FMM
\cite{ambikasaran2014inverse,coulier2017inverse}, orthogonal-basis
factorization \cite{ma2017accuracy,Ma2019,boukaram2025linear}, and randomized
factorization \cite{yesypenko2025}, interleaves low-rank compression
with elimination: a box retains a small set of skeleton degrees of freedom,
eliminates its redundant degrees of freedom, and passes its skeletons to
the next coarser level. This approach can attain
linear complexity when numerical ranks remain bounded. The spatial hierarchy exposes parallelism among boxes, but exploiting it in distributed memory requires a schedule that respects elimination dependencies and controls communication costs.

Existing parallelization efforts include shared-memory implementations \cite{takahashi2020parallelization,boukaram2025linear}, $\mathcal{H}^2$-ULV methods that incorporate fill-ins into shared bases to expose parallelism \cite{ma2022scalable,Ma_2024}, and iterative hierarchical solvers \cite{spendlhofer2025interplay,sushnikova2026mg}. The closest distributed-memory RS-S implementation is the 2D solver of Liang et al.\ \cite{liang2024on}. It partitions the 2D domain among processes and separates interior from boundary elimination. A four-color process-grid schedule lets approximately one quarter of the processes eliminate boundary boxes during each stage.

\begin{table}[t]
\caption{Per-box stencil and scheduling quantities in 2D and 3D.}
\label{t:2d-vs-3d}
\centering
\footnotesize
\begin{tabular}{@{}lrr@{}}
\toprule
Quantity & 2D & 3D\\
\midrule
Near-field boxes $|\N|$ & 8 & 26\\
Distance-2 boxes $|\M|$ & 16 & 98\\
Off-diagonal neighbor pairs receiving updates & 28 & 325\\
Children per refinement & 4 & 8\\
Active processes per stage, process coloring & 1 in 4 & 1 in 8\\
Boundary share, $v=8$ / $v=4$ boxes per side & 44\% / 75\% & 58\% / 88\%\\
\bottomrule
\end{tabular}
\end{table}

Three-dimensional problems make storage, communication, and scheduling substantially harder than in 2D (\Cref{t:2d-vs-3d}). First, each box has $26$ near-field and $98$ distance-2 boxes instead of $8$ and $16$. An eliminated box therefore contributes Schur updates to $325$ off-diagonal neighbor pairs instead of $28$, and its compression reads a much larger target, which drives temporary storage and communication volume. Second, each refinement multiplies the number of boxes by eight, so keeping each process's subdomain large enough to limit its boundary share requires using many cores per process. Third, numerical rank growth increases coarse-level costs: the average skeleton size in our 3D Laplace runs at $\epsilon=10^{-3}$ roughly doubles from the leaf level to three levels up (\Cref{t:exp-level-breakdown-omp128}). The 2D ranks reported in~\cite[Fig.~9]{liang2024on}, using the tighter default tolerance $\epsilon=10^{-6}$, are smaller and grow more slowly. Different tolerances and experimental settings prevent isolating the effect of dimensionality. Coarse levels thus have few boxes with large active sets, which limits concurrency where boundary boxes dominate each subdomain. Finally, a direct 3D extension of process-grid coloring needs eight colors, so only one in eight processes would eliminate boundary boxes during each stage, causing even more severe process idling issues than 2D. 

We present a spatial-coloring-based distributed $\mathcal{H}^2$ factorization that addresses these challenges, with the following contributions.
\begin{enumerate}
\item \textit{Spatial coloring for distributed elimination.}
We color the spatial boxes, rather than the process grid, by the parity of their grid coordinates. Each process owns boxes of every color and can contribute eligible boundary boxes during every color stage. The same colors schedule shared-memory parallelism within each process.

\item \textit{Storage and communication reduction for 3D stencils.}
We reduce the storage and communication costs of the Schur complement updates by delaying the GEMM products until they are needed by the process that owns the target block. We also reduce the storage costs of the interpolative decomposition by streaming the sketching operation one block at a time, instead of assembling the full compression target.

\item \textit{Hybrid parallelism for box growth and coarse-level ranks.}
Shared-memory parallelism over color batches lets one process per node use all of its cores at fine levels. At coarse levels, cooperative tasks within a box, right-side BLAS-3 solves, and process reduction keep the cores busy despite few boxes and large ranks.

\item \textit{Large-scale evaluation.}
We demonstrate factorization of a 3D Laplace system with $N=1536^3=3,623,878,656$ unknowns on 512 nodes, alongside 3D Helmholtz tests and a 2D Mat\'ern system with $N=16384^2=268,435,456$ unknowns. To our knowledge, this is the first distributed-memory strong recursive skeletonization factorization demonstrated for three-dimensional volume integral equations at billion-unknown scale.
\end{enumerate}

\Cref{sec:background,sec:dist-algorithms} introduce RS-S and the distributed schedule. \Cref{sec:lazy-schur,sec:streaming-id,sec:shared-memory} present the optimizations; \Cref{sec:alg-analysis,s:experiments} analyze costs and report experiments.

%% file: background.tex
\section{Background}
\label{sec:background}

This section reviews the strong recursive skeletonization (RS-S) factorization~\cite{minden2017recursive} and its data dependencies in the distributed-memory setting~\cite{liang2024on}, and fixes the notation used in later sections. We state only the operations that our algorithm and optimizations act on and refer to~\cite{minden2017recursive,liang2024on} for their derivations; related factorizations include~\cite{ambikasaran2014inverse,coulier2017inverse,pouransari2017fast,sushnikova2018compress,sushnikova2023fmm}.

Throughout, we consider a point set $\X = \{\bm{x}_i\}_{i=1}^N$ in $\mathbb{R}^d$ ($d=2$ or $3$) and a symmetric kernel $K(\bm{x},\bm{y})=K(\bm{y},\bm{x})$ defining $A_{i,j}=K(\bm{x}_i,\bm{x}_j)$ for $i\ne j$. We assume $A=A^\top$ throughout the factorization, storage model, and parallel algorithms. Examples include real-symmetric Laplace and Mat\'ern matrices and the complex-symmetric Helmholtz matrix. The superscript $\top$ always denotes the non-conjugating transpose, and $A_{\I,\J}$ denotes the submatrix of $A$ with row indices $\I$ and column indices $\J$.

%=============================================================================
\subsection{Hierarchical Domain Decomposition}
\label{subsec:hierarchy}
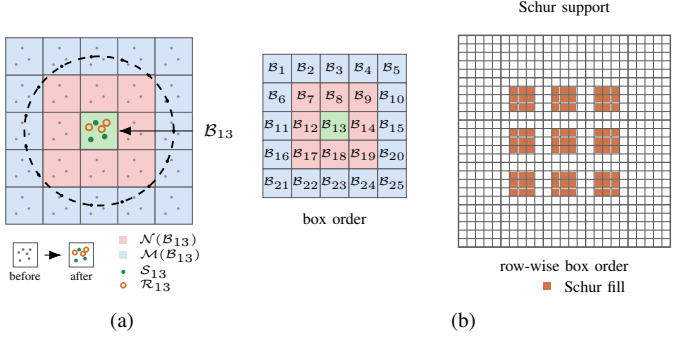
\begin{figure}[t]
     \centering
     \subfloat[]{\begin{minipage}[b]{0.36\linewidth}
         \centering
         \resizebox{\linewidth}{!}{%
         \begin{tikzpicture}[scale=0.54, every node/.style={font=\scriptsize}]
             \definecolor{centerbox}{RGB}{200,232,194}
             \definecolor{nearbox}{RGB}{247,204,204}
             \definecolor{distbox}{RGB}{214,228,248}
             \definecolor{gridline}{RGB}{92,92,92}
             \definecolor{skelgreen}{RGB}{31,136,78}
             \definecolor{redundantorange}{RGB}{213,113,38}
             \foreach \i in {0,...,4} {
                 \foreach \j in {0,...,4} {
                     \fill[distbox] (\i,\j) rectangle +(1,1);
                     \draw[gridline, line width=0.35pt] (\i,\j) rectangle +(1,1);
                     \fill[black!45] (\i+0.28,\j+0.28) circle (0.035);
                     \fill[black!45] (\i+0.62,\j+0.42) circle (0.035);
                     \fill[black!45] (\i+0.44,\j+0.72) circle (0.035);
                 }
             }
             \foreach \i in {1,...,3} {
                 \foreach \j in {1,...,3} {
                     \fill[nearbox] (\i,\j) rectangle +(1,1);
                     \draw[gridline, line width=0.35pt] (\i,\j) rectangle +(1,1);
                     \fill[black!45] (\i+0.28,\j+0.28) circle (0.035);
                     \fill[black!45] (\i+0.62,\j+0.42) circle (0.035);
                     \fill[black!45] (\i+0.44,\j+0.72) circle (0.035);
                 }
             }
             \fill[centerbox] (2,2) rectangle +(1,1);
             \draw[gridline, line width=0.55pt] (2,2) rectangle +(1,1);
             \draw[black, line width=0.75pt, dashed] (2.5,2.5) circle (2.0);
             \foreach \a in {0,24,...,336} {
                 \fill[black] ($(2.5,2.5)+(\a:2.0)$) circle (0.045);
             }
             \fill[skelgreen] (2.28,2.27) circle (0.075);
             \fill[skelgreen] (2.65,2.35) circle (0.075);
             \fill[skelgreen] (2.41,2.72) circle (0.075);
             \draw[redundantorange, line width=0.75pt] (2.70,2.72) circle (0.085);
             \draw[redundantorange, line width=0.75pt] (2.23,2.60) circle (0.085);
             \draw[redundantorange, line width=0.75pt] (2.57,2.55) circle (0.085);
             \node[font=\scriptsize, anchor=west] at (5.08,2.50) {$\B_{13}$};
             \draw[-{Latex[length=1.8mm]}, line width=0.4pt] (5.00,2.50) -- (3.02,2.50);

             \draw[gridline, line width=0.35pt] (0.20,-1.15) rectangle +(0.65,0.65);
             \foreach \p in {(0.382,-0.975),(0.623,-0.923),(0.467,-0.682),(0.655,-0.682),(0.350,-0.760),(0.571,-0.793)} {
                 \fill[black!55] \p circle (0.035);
             }
             \node[font=\tiny] at (0.52,-1.40) {before};
             \draw[-{Latex[length=1.8mm]}, line width=0.4pt] (1.04,-0.83) -- (1.50,-0.83);
             \draw[gridline, line width=0.35pt] (1.70,-1.15) rectangle +(0.65,0.65);
             \fill[skelgreen] (1.882,-0.975) circle (0.05);
             \fill[skelgreen] (2.123,-0.923) circle (0.05);
             \fill[skelgreen] (1.967,-0.682) circle (0.05);
             \draw[redundantorange, line width=0.65pt] (2.155,-0.682) circle (0.06);
             \draw[redundantorange, line width=0.65pt] (1.850,-0.760) circle (0.06);
             \draw[redundantorange, line width=0.65pt] (2.071,-0.793) circle (0.06);
             \node[font=\tiny] at (2.03,-1.40) {after};

             \fill[nearbox] (3.05,-0.52) rectangle +(0.20,0.20);
             \node[font=\tiny, anchor=west] at (3.36,-0.42) {$\N(\B_{13})$};
             \fill[distbox] (3.05,-0.94) rectangle +(0.20,0.20);
             \node[font=\tiny, anchor=west] at (3.36,-0.84) {$\M(\B_{13})$};
             \fill[skelgreen] (3.15,-1.28) circle (0.055);
             \node[font=\tiny, anchor=west] at (3.36,-1.28) {$\S_{13}$};
             \draw[redundantorange, line width=0.65pt] (3.15,-1.64) circle (0.065);
             \node[font=\tiny, anchor=west] at (3.36,-1.64) {$\R_{13}$};
         \end{tikzpicture}%
         }
         
     \end{minipage}}%
     \hfill%
     \subfloat[]{\begin{minipage}[b]{0.62\linewidth}
         \centering
         \resizebox{\linewidth}{!}{%
         \begin{tikzpicture}[every node/.style={font=\scriptsize}]
             \definecolor{centerbox}{RGB}{200,232,194}
             \definecolor{nearbox}{RGB}{247,204,204}
             \definecolor{distbox}{RGB}{214,228,248}
             \definecolor{gridline}{RGB}{92,92,92}
             \definecolor{fillorange}{RGB}{201,91,42}
             \def\b{0.44}
             \def\mx{3.00}
             \def\my{0.02}
             \def\ms{0.13}
             \foreach \i/\j/\lab/\col in {
                 0/4/1/distbox,1/4/2/distbox,2/4/3/distbox,3/4/4/distbox,4/4/5/distbox,
                 0/3/6/distbox,1/3/7/nearbox,2/3/8/nearbox,3/3/9/nearbox,4/3/10/distbox,
                 0/2/11/distbox,1/2/12/nearbox,2/2/13/centerbox,3/2/14/nearbox,4/2/15/distbox,
                 0/1/16/distbox,1/1/17/nearbox,2/1/18/nearbox,3/1/19/nearbox,4/1/20/distbox,
                 0/0/21/distbox,1/0/22/distbox,2/0/23/distbox,3/0/24/distbox,4/0/25/distbox} {
                 \fill[\col] (\i*\b,\j*\b+0.78) rectangle +(\b,\b);
                 \draw[gridline, line width=0.3pt] (\i*\b,\j*\b+0.78) rectangle +(\b,\b);
                 \node[font=\tiny] at ({\i*\b+0.5*\b},{\j*\b+0.78+0.5*\b}) {$\B_{\lab}$};
             }
             \node at (1.10,0.45) {box order};

             \foreach \ii/\xi/\yi in {
                 7/1/3,8/2/3,9/3/3,12/1/2,13/2/2,14/3/2,17/1/1,18/2/1,19/3/1} {
                 \foreach \jj/\xj/\yj in {
                     7/1/3,8/2/3,9/3/3,12/1/2,13/2/2,14/3/2,17/1/1,18/2/1,19/3/1} {
                     \fill[fillorange!85]
                         ({\mx+(\jj-1)*\ms+0.018},{\my+(25-\ii)*\ms+0.018})
                         rectangle +(0.094,0.094);
                 }
             }
             \foreach \k in {0,...,25} {
                 \draw[gridline, line width=0.17pt]
                     ({\mx+\k*\ms},{\my}) -- ({\mx+\k*\ms},{\my+25*\ms});
                 \draw[gridline, line width=0.17pt]
                     ({\mx},{\my+\k*\ms}) -- ({\mx+25*\ms},{\my+\k*\ms});
             }
             \draw[gridline, line width=0.6pt] (\mx,\my) rectangle +(3.25,3.25);
             \foreach \k in {5,10,15,20} {
                 \draw[gridline, line width=0.45pt]
                     ({\mx+\k*\ms},{\my}) -- ({\mx+\k*\ms},{\my+25*\ms});
                 \draw[gridline, line width=0.45pt]
                     ({\mx},{\my+\k*\ms}) -- ({\mx+25*\ms},{\my+\k*\ms});
             }
             \node at ({\mx+12.5*\ms},{\my-0.28}) {row-wise box order};
             \node at ({\mx+12.5*\ms},{\my+25*\ms+0.42}) {Schur support};

             \fill[fillorange!85] ({\mx+1.30},{\my-0.66}) rectangle +(0.12,0.12);
             \node[anchor=west] at ({\mx+1.50},{\my-0.60}) {Schur fill};
         \end{tikzpicture}%
         }
         
     \end{minipage}}
     \caption{Proxy compression and local elimination dependencies. (a) Compression of $\B_{13}$ uses distance-2 information, proxy points (on the circle), and an ID that splits points in $\B_{13}$ into skeleton indices $\S_{13}$ and redundant indices $\R_{13}$. (b) For readability, the 25-box stencil is numbered row-wise; the implementation uses Morton ordering. Assuming $\B_{13}$ is eliminated first, eliminating $\R_{13}$ creates Schur-complement contributions on $(\S_{13}\cup\N(\B_{13}))\times(\S_{13}\cup\N(\B_{13}))$ as in~\eqref{eq:schur}. In the box-level matrix, the $\B_{13}$ row and column represent the remaining skeleton variables $\S_{13}$; colored entries mark the blocks receiving this contribution.  
     }
     \label{f:proxy}
\end{figure}

%=============================================================================
Assume the points lie in a $d$-dimensional unit cube (a unit square in 2D). We recursively subdivide each box into $2^d$ equal-sized children, forming a $2^d$-ary tree $\T$ (a quadtree in 2D or an octree in 3D). Our focus is three-dimensional point sets, whose octree boxes have substantially larger near-field and distance-2 stencils than in 2D, as quantified below. The method also applies to $d=2$, and \Cref{f:proxy,f:color-waves} use 2D examples for readability. For the algorithm and experiments considered here, we assume a uniform point distribution and use a \emph{perfect} tree: the root is at level $0$ and all leaves are at level $L$, with a bounded number $n_{\mathrm{leaf}}=\bigO(1)$ of points per leaf. The coordinate-parity coloring schedule described in \Cref{sec:dist-algorithms} is restricted to this uniform-tree setting.

We refer to each node $\B$ of $\T$ as a \emph{box}. Let $n_\ell = v_\ell^d$ denote the number of boxes at level $\ell$. A box centered at $((i_1+\frac{1}{2})\Delta,\ldots,(i_d+\frac{1}{2})\Delta)$ with $0\leq i_k<v_\ell$ and $\Delta=1/v_\ell$ has integer coordinates $(i_1,\ldots,i_d)$ and a global Morton key $k_m(\B)$ obtained by interleaving the binary bits of $i_1,\ldots,i_d$. For two boxes at the same level, $\operatorname{dist}_{\infty}(\B,\B')=\max_{j}|i_j(\B)-i_j(\B')|$ measures their distance in units of the box width.

We denote the $i$th box at level $\ell$ by $\B^\ell_i$ and drop the superscript and subscript when clear from context. The \emph{parent} $\P(\B)$ and \emph{children} $\C(\B)$ carry their usual tree meanings. The \emph{near-field} $\N(\B)$ consists of the same-level boxes with $\operatorname{dist}_{\infty}(\B,\B')=1$, namely those sharing a face, edge, or corner with $\B$; there are at most $3^d-1$ of them ($8$ in 2D, $26$ in 3D). The \emph{far-field} $\F(\B)$ consists of all remaining same-level boxes, and the \emph{distance-2 far-field} $\M(\B)\subset\F(\B)$ consists of those with $\operatorname{dist}_{\infty}(\B,\B')=2$; there are at most $5^d-3^d$ of them ($16$ in 2D, $98$ in 3D). With a slight abuse of notation, we also use $\B$, $\N$, $\F$, and $\M$ to denote the point indices contained in the corresponding boxes. \Cref{f:proxy}(a) illustrates these sets.

%=============================================================================
\subsection{Strong Skeletonization of a Box}
\label{subsec:lowrank}\label{subsec:elimination}
%=============================================================================
For suitable smooth kernels and fixed separation geometry, interactions between well-separated point sets have numerical ranks that remain bounded independently of $N$ at a prescribed accuracy $\eps$~\cite{minden2017recursive}. RS-S exploits this structure through the column interpolative decomposition (ID): given $A_{\I,\J}\in\mathbb{C}^{m\times n}$, it selects \emph{skeleton} indices $\S\subset\J$, \emph{redundant} indices $\R=\J\setminus\S$, and an interpolation matrix $T\in\mathbb{C}^{|\S|\times|\R|}$ such that
\begin{equation}
    \|A_{\I, \R} - A_{\I, \S} \, T\|
     \le \eps \, \|A_{\I, \J}\|.\label{eq:id-approx}
\end{equation}
A rank-revealing QR (RRQR) factorization computes the ID at $\bigO(mnk)$ cost with rank $k = |\S|$~\cite{gu1996efficient,cheng2005compression}.

\paragraph{Compression}
Eliminating earlier boxes modifies the interactions between $\B$ and $\M(\B)$ through Schur complement updates, whereas the interactions with $\F(\B)\setminus\M(\B)$ remain pure kernel evaluations. Compression therefore applies the ID to
\begin{equation}
\label{eq:proxy}
W_\B=\begin{pmatrix}
A_{\M(\B),\B}\\
K_{\mathrm{proxy},\B}
\end{pmatrix},
\end{equation}
where $A_{\M(\B),\B}$ is the current, possibly updated block and $K_{\mathrm{proxy},\B}$ contains kernel evaluations between $\B$ and a small set of points on a \emph{proxy surface} surrounding $\B$, which captures the column space of the remaining far-field interactions~\cite{minden2017recursive}. The validity of this representation depends on the kernel; \Cref{s:exp-acc} discusses the limitation encountered for the Mat\'ern kernel. \Cref{sec:streaming-id} describes how we sketch $W_\B$ without materializing it in full.

\paragraph{Elimination}
Given the ID of $\B$ with skeleton $\S$, redundant indices $\R=\B\setminus\S$, and interpolation matrix $T$, a unit-triangular transformation defined by $T$ decouples $\R$ from the far field up to the ID accuracy~\cite{minden2017recursive}. Only three distinct transformed blocks need to be formed:
\begin{equation}
\label{eq:Xblocks}
\begin{aligned}
X_{\R,\R}={}&A_{\R,\R}-T^{\top}A_{\S,\R}-A_{\S,\R}^{\top}T\\
            &+T^{\top}A_{\S,\S}T,\\
X_{\S,\R}={}&A_{\S,\R}-A_{\S,\S}T,\\
X_{\N,\R}={}&A_{\N,\R}-A_{\N,\S}T,
\end{aligned}
\end{equation}
with $X_{\R,\R}=X_{\R,\R}^{\top}$, $X_{\R,\S}=X_{\S,\R}^{\top}$, and $X_{\R,\N}=X_{\N,\R}^{\top}$. After this \emph{sparsification} step, the redundant indices interact only with $\S\cup\N$. Let $D_\B=X_{\R,\R}$ be the symmetric redundant pivot block and assume it is nonsingular. Symmetric block elimination of $\R$ then produces the Schur complement updates
\begin{equation}
\label{eq:schur}
\begin{aligned}
\tilde A_{\S,\S}&=A_{\S,\S}-X_{\S,\R}D_\B^{-1}X_{\S,\R}^{\top},\\
\tilde A_{\N,\S}&=A_{\N,\S}-X_{\N,\R}D_\B^{-1}X_{\S,\R}^{\top},\\
\tilde A_{\N,\N}&=A_{\N,\N}-X_{\N,\R}D_\B^{-1}X_{\N,\R}^{\top},
\end{aligned}
\end{equation}
and leaves the $(\S,\F)$, $(\N,\F)$, and $(\F,\F)$ blocks unchanged. We factor $D_\B$ with a pivoted factorization and apply $D_\B^{-1}$ through solves with the stored factors (\Cref{subsec:blas3}); for complex-symmetric blocks, the factorization uses a non-conjugating transpose. The per-box cost is $\bigO(k^3)$, provided the redundant and active neighbor sizes are also $\bigO(k)$.

We denote by $V_\B$ the product of the permutation that splits $\B$ into $\R$ and $\S$, the interpolation transformation, and the block-elimination transformation, so that, in the compressed approximation, $V_\B A V_\B^{\top}$ is block diagonal between $\R$, with pivot block $D_\B$, and the remaining indices, which carry the updated blocks~\eqref{eq:schur}. Since $V_\B$ consists of a permutation and unit-triangular transformations, its inverse and transpose can be applied from the same stored data. The inverses of the elementary unit-triangular transformations are applied by negating their off-diagonal blocks and reversing the transformation order, without forming dense inverse matrices.

Let $A^\ell_i$ denote the active matrix just before processing $\B^\ell_i$, so that $A^\ell_{i+1}=V_{\B^\ell_i}A^\ell_iV_{\B^\ell_i}^{\top}$. After all $n_\ell$ boxes at level $\ell$ have been processed, the remaining active block $\tilde A^\ell$ is supported on the surviving skeleton indices $\bigcup_{i=1}^{n_\ell}\S_i$.

%=============================================================================
\subsection{Multi-Level Factorization}
\label{subsec:multilevel}
%=============================================================================
At level $\ell-1$, each parent box owns the skeleton indices of its children, and compression and elimination are repeated at levels $L,L-1,\ldots,1$ (\Cref{a:fact}). The residual $\tilde A^1$ is assembled at the root and factored by a dense pivoted factorization. Let $\mathcal V=\prod_{\ell=1}^{L}\prod_{i=n_\ell}^{1}V_{\B_i^\ell}$ be the product of all box transformations in application order, and let $D$ be the block-diagonal matrix of the pivot blocks $D_{\B_i^\ell}$ and the root residual $\tilde A^1$. The approximate factorization and its inverse are
\begin{equation}
\label{eq:approx-factor}
F\equiv\mathcal V^{-1}D\mathcal V^{-\top}\approx A,
\end{equation}
\begin{equation}
\label{eq:approx-inv}
F^{-1}=\mathcal V^{\top}D^{-1}\mathcal V\approx A^{-1}.
\end{equation}
Applying $F^{-1}$ to a vector uses three stages: an upward pass applying the $V_{\B_i^\ell}$ transformations in elimination order, a block-diagonal solve with the redundant pivot factors and the dense root factors, and a downward pass applying $V_{\B_i^\ell}^{\top}$ in reverse elimination order.

Strong recursive skeletonization also applies to nonsymmetric matrices. In that setting, compression must account for both row and column interactions, and elimination uses distinct left and right transformations. The formulation and implementation presented here assume symmetry and exploit it in compression, factor storage, and update computation.

\input{factor.tex}

%=============================================================================
\subsection{Dependency Structure of $\mathcal{H}^2$ Factorization}
\label{subsec:dependency}\label{sec:data_dep}
%=============================================================================
The data read and written by one call to \textsc{StrongSkeletonize}$(\B)$ were analyzed in~\cite{liang2024on}; we restate the facts used in later sections. Compression reads the current distance-2 blocks $A_{\M(\B),\B}$ and the proxy interactions in~\eqref{eq:proxy}. Unchanged kernel blocks can be generated locally from point coordinates, but a block $A_{C,\B}$ modified by the earlier elimination of a common neighbor $E$ of $C$ and $\B$ must also include the contribution of $E$, supplied either as an explicit block or as a generator (\Cref{sec:lazy-schur}). Elimination reads $A_{\B,\B}$ and $A_{\N(\B),\B}$ and writes the Schur updates~\eqref{eq:schur} on $(\S\cup\N(\B))\times(\S\cup\N(\B))$, whose box pairs are at most distance two apart. Hence a box must receive the contributions of previously eliminated neighbors before it is compressed or eliminated, and an elimination never modifies interactions beyond distance two.

%% file: factor.tex
\begin{algorithm}[t]
\caption{Sequential strong-skeletonization factorization}
\label{a:fact}
\begin{algorithmic}[1]
\Require{Symmetric kernel $K$, points $\X = \{x_i\}_{i=1}^N \in \mathbb{R}^{d}$.}
\Ensure Factors in~\eqref{eq:approx-factor}.
\State Build the spatial tree of \Cref{subsec:hierarchy}, with root level $0$ and leaf level $L$, and assign points to leaf boxes.
\For{$\ell = L, L-1, \ldots, 1$}
    \For{$i = 1, 2, \ldots, n_\ell$}\label{l:inloop}
        \State \Call{StrongSkeletonize}{$\B_i^\ell$}
    \EndFor
    \State Assemble the next-level active matrix from surviving child skeletons.
\EndFor
\State Factor the dense residual $\tilde A^1$ assembled at the root and store its pivot factors.

\Function{StrongSkeletonize}{$\B_i^\ell$}
    \State Compute an ID of the proxy compression target formed from the active matrix $A^\ell_i$, giving $\S$, $\R$, and $T$ as in~\eqref{eq:id-approx}.
    \State Form the ID-decoupled blocks $X_{\R,\R}$, $X_{\S,\R}$, $X_{\N,\R}$ of~\eqref{eq:Xblocks}.
    \State Factor the symmetric pivot block $D_{\B^\ell_i}=X_{\R,\R}$.
    \State Apply the Schur complement updates~\eqref{eq:schur} to eliminate $\R$ and update the active blocks on $\S \cup \N$.
    \State Store $V_{\B^\ell_i}$ and the pivot factors of $D_{\B^\ell_i}=X_{\R,\R}$.
\EndFunction
\end{algorithmic}
\end{algorithm}

% The data dependency of one call to \textsc{StrongSkeletonize}$(\B)$ is what the
% distributed algorithms of \Cref{sec:dist-algorithms} are built around, and we
% record it here for reference: the compression reads the existing distance-2
% interactions $A_{\M,\B}$; the elimination reads the active
% near-field blocks $A_{\B,\B}$ and $A_{\N,\B}$, with opposite orientations obtained by transpose; and the resulting
% Schur complement updates among other boxes modify blocks $A_{\N,\N}$ only. \Cref{sec:data_dep}
% develops this into the dependency structure used later.

%% file: distributed-algorithm.tex
\section{Distributed Algorithm}
\label{sec:dist-algorithms}

We propose a distributed-memory $\mathcal{H}^2$ factorization algorithm that uses coordinate-parity coloring of spatial boxes to schedule boundary elimination, referred to as the spatial-coloring algorithm. Unlike process-grid coloring \cite{liang2024on}, which in 3D would leave seven of every eight processes idle during each boundary stage, it lets every process with eligible boundary boxes participate in each stage simultaneously. It runs the multilevel sweep of Algorithm~\ref{a:fact} from the leaf level to the root level. At each level, it uses a regular block partition of the domain. Specifically, for $n_l=v_l^d$ boxes at level $l$, each process $p$ owns a subdomain containing contiguous $v\times v$ boxes in 2D or $v\times v \times v$ boxes in 3D, %
denoted by $S_p$ (see \cref{f:color-waves}(b) for an example of 4 processes in 2D).

At a single level, a boundary box is a box in $S_p$ adjacent to a box owned by another process; let $\Gamma_p$ be the set of boundary boxes. The remaining boxes $S_p\setminus\Gamma_p$ are interior boxes. Boundary-box compression and elimination can read modified near-field and distance-2 interactions across process boundaries. An interior box has only local near-field neighbors, although its distance-2 compression interactions can still involve remote boxes.

The algorithm partitions boundary boxes into $2^d$ \emph{color sets} and processes them in color order. Boxes of one color are nonadjacent and can be processed concurrently, with shared-target updates handled separately as discussed in \Cref{subsec:dependency,sec:shared-memory}. After each boundary color, processes exchange the update representations needed by subsequent consumers. Interior boxes are processed after the final boundary exchange, under the data-availability requirement described below.

\subsubsection{Graph Coloring}
\begin{figure}[t]
\centering
\includegraphics[width=\linewidth]{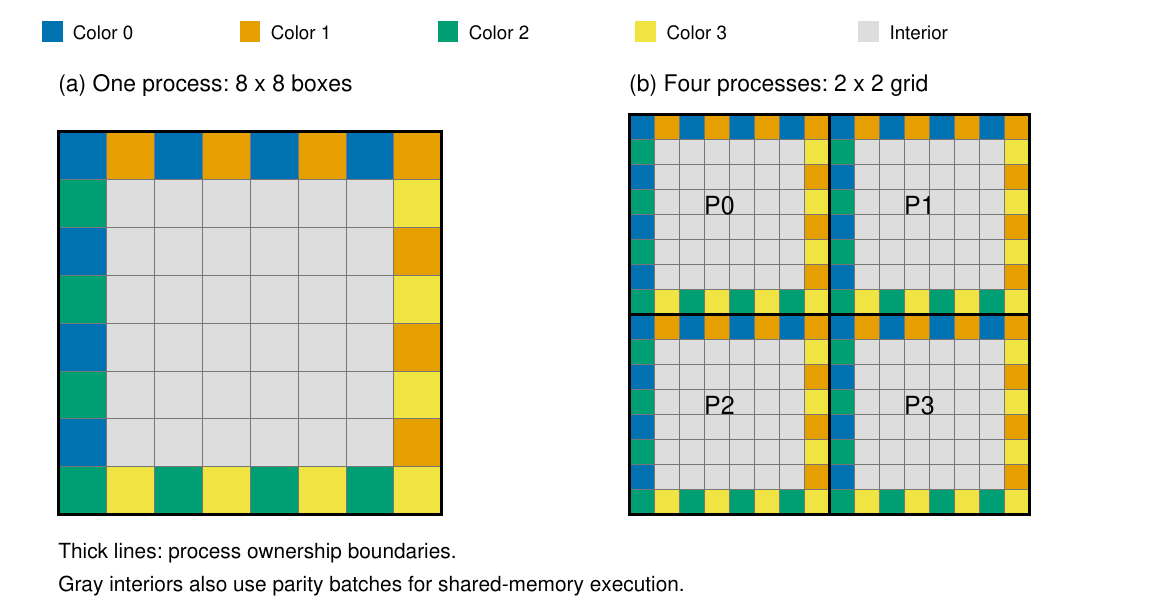}
\caption{Coordinate-parity coloring for the spatial-coloring algorithm ($v=8$, $d=2$). \textbf{(a)}~One process owns an $8\times8$ subdomain. Boundary boxes show their parity colors; gray interior boxes also use parity batches for shared-memory execution. \textbf{(b)}~A $2\times2$ process grid. Thick lines indicate process ownership boundaries. Colors use global box coordinates and remain consistent across processes.}
\label{f:color-waves}
\end{figure}

Using the integer box coordinates of \Cref{subsec:hierarchy}, assign each box the color
\[
w(\B)=\sum_{j=1}^{d}2^{j-1}(i_j\bmod 2)
\in\{0,1,\ldots,2^d-1\}.
\]
This gives four colors in 2D and eight in 3D. Distinct same-color boxes agree in every coordinate parity, so their grid distance is at least two; they are therefore nonadjacent. By the dependencies in \Cref{subsec:dependency}, eliminating one does not change the blocks read to compress or eliminate another box of the same color. Their update supports can nevertheless overlap on shared neighbor pairs. Concurrent processing therefore requires separate handling of these contributions before a later box consumes the target block; \Cref{sec:shared-memory} describes deferred accumulation by target pair.
The color can also be obtained as $w(\B)=k_m(\B)\bmod 2^d$ from the global Morton key. \Cref{f:color-waves} illustrates the boundary colors in 2D.

% We use the same parity classes for shared-memory scheduling within each
% process (\Cref{sec:shared-memory}). Their low-cost Morton-key construction
% therefore serves both MPI and OpenMP parallelism. Other proper colorings are
% possible, but finding fewer colors on a process boundary is not our objective;
% we target the process concurrency and data-movement costs of a regular,
% locally computable schedule.

\subsubsection{Elimination Order}
The spatial-coloring algorithm uses the same workflow as the sequential $\mathcal{H}^2$ factorization algorithm~\ref{a:fact}, except that the per-level box loop at Line \ref{l:inloop} is replaced by Algorithm~\ref{a:color}.

At level $l$, boundary boxes $\Gamma_p$ on each process $p$ are eliminated in the order
\begin{eqnarray}
	\label{eq:elim-order}
	e(w=0) <e(w=1)<\cdots< e(w=2^d-1),\label{e:oc_c}
\end{eqnarray}
where a smaller value of $e(\cdot)$ means earlier elimination. Boxes within a color can be processed concurrently with deferred accumulation of shared-target updates to avoid race condition (\Cref{sec:shared-memory}). After the final boundary exchange, each process eliminates its interior boxes. We also use the color order~\eqref{e:oc_c} for interior boxes to expose shared-memory parallelism (\Cref{sec:shared-memory}).

\begin{algorithm*}[t]
\caption{Spatial-coloring algorithm: distributed elimination, run at levels $L$ through $1$ before the dense root factorization. Tag colors: \textcolor{commBdy}{boundary communication}, \textcolor{commLvl}{level-transition communication}.}
\label{a:color}
\begin{algorithmic}[1]
\Require Boxes at the current level $\ell$ under the block partition with $P$ processes, symmetric kernel $K$, tolerance $\eps$.
\Ensure All boxes at this level eliminated and $A^{\ell-1}_1$ at next level assembled.
\State assign each local box $\B$ to color $w(\B) = k_m(\B) \bmod 2^d$
\For{$w = 0, 1, \dots, 2^d-1$} \algcmt{color sets for the boundary boxes}\label{l:b_color}
    \State eliminate \emph{boundary} boxes $\B$ with $w(\B) = w$ by \Call{StrongSkeletonize}{$\B$}\label{l:s_c_1}
    \State \textcolor{commBdy}{\textbf{exchange}} the resulting update representations with the $q\le3^d-1$ neighboring processes; install them before the next color or the interior phase\label{l:b_com}
\EndFor
\For{$w = 0, 1, \dots, 2^d-1$} \algcmt{color batches for the interior boxes}\label{l:i_color}
    \State eliminate \emph{interior} boxes $\B$ with $w(\B) = w$ by \Call{StrongSkeletonize}{$\B$}\label{l:s_c_2}
\EndFor
\State assemble $A_1^{\ell-1}$ from surviving child skeletons, consuming pending generators
\State \textcolor{commLvl}{\textbf{(only if $n_{\ell-1}/P<n_{\min}$) redistribute}} $A^{\ell-1}_1$ with a reduced process count
\end{algorithmic}
\end{algorithm*}

After all boxes at level $\ell$ have been processed, surviving child skeletons are assembled into $A^{\ell-1}_1$. If $n_{\ell-1}/P<n_{\min}$, the algorithm reduces the active process count and redistributes the next-level matrix so that each active process owns sufficiently many boxes. At level $0$, the dense root matrix is factored as in Algorithm~\ref{a:fact}.

\subsubsection{Communication Pattern}
Before a dependent boundary color starts, all required updates from earlier colors must have been received and processed. Line~\ref{l:b_com} performs an exchange after each of the $2^d$ boundary colors, including the last color. Before interior processing begins, the final boundary exchange completes delivery of all remote distance-2 blocks or update representations needed by interior compression, including contributions from earlier colors. Their coordinates and active-index information are available locally, so the interior phase requires no further communication. Interior elimination then has local near-field support, so its newly generated updates do not cross process boundaries.

Elimination of one boundary box $\B$ produces Schur complement updates to each of the $\bigo{1}$ box pairs in $(\N,\N)$, where a pair $(A,C)$ requires $\bigo{b_A b_C}$ scalar words for its Schur update, with $b_A,b_C=\bigo{k}$ the active box sizes, which may or may not have been skeletonized/sparsified. In our baseline eager implementation, the dense Schur updates are computed after each color and buffered for exchange with $q \le 3^d - 1$ neighboring processes, which can lead to large communication volume and high peak memory usage. This is particularly problematic in 3D, where the number of off-diagonal neighbor pairs is $325=\bigo{26^2/2}$ (even accounting for symmetry) each requiring $\bigo{k^2}$ temporary storage units. Therefore, we propose a lazy (or deferred) Schur update strategy in Section~\ref{sec:lazy-schur} to reduce the communication volume and peak memory usage.

%% file: lazy-schur.tex
\section{On-Demand Schur Updates}
\label{sec:lazy-schur}
We reduce the storage and communication costs of the Schur complement updates. These costs grow sharply in 3D (\Cref{t:2d-vs-3d}): an eliminated box updates up to $325$ off-diagonal neighbor pairs in 3D, compared with $28$ in 2D.

\subsection{A generator representation for Schur updates}
\label{subsec:generators}
For an eliminated source box $E$, let
$b_E=|\S_E|+|\R_E|$, $k_E=|\S_E|$, and $r_E=|\R_E|$. Let $D_E=X_{\R_E,\R_E}$ denote the
unfactored redundant diagonal block from~\eqref{eq:Xblocks}. Define the \textit{generator}
\begin{equation}
\label{eq:schur-generator}
\begin{aligned}
 G_E&=-X_{\N(E),\R_E}D_E^{-1},\\
 G_{E,A}&=G_E[A,:],\quad A\in\N(E).
\end{aligned}
\end{equation}
%\chao{Does notation $ G_{E,A}$ make more sense above?}
The Schur update of $E$ to
box pair $A,C\in\N(E)$ is
\begin{equation}
\label{eq:generated-pair}
 \Delta^{(E)}_{A,C}
 =-X_{A,\R_E}D_E^{-1}X_{C,\R_E}^{\top}
 =-G_{E,A}D_EG_{E,C}^{\top}.
\end{equation}
The retained source representation includes $D_E$, the
required generator $G_E$, and skeleton-coupling data needed when either $A$ or $C$ has already been skeletonized. The eager implementation communicates $\Delta^{(E)}_{A,C}$. In contrast, the on-demand implementation communicates $D_E$ and $G_{E}$ to the recipient of the Schur update. The products are formed on demand for three purposes: later box elimination, ID compression, and parent assembly.

\subsection{Two evaluation policies}
\label{subsec:lazy-modes}
We use a parameter \texttt{lazy\_schur} to select when the Schur updates are evaluated. \texttt{lazy\_schur}=0 and 1 correspond to the two evaluation modes below.

\paragraph{Eager update (mode 0)}
The source process forms all dense products as soon as a
color set completes. All products required by remote processes are buffered and exchanged before the next dependent color.
Compression and parent construction read the stored modified blocks. This
policy avoids later regeneration but stores both near-field and distance-2 far-field blocks and can retain many products in pending send buffers.

\paragraph{On-demand near-field and far-field update (mode 1)}
Both near-field and distance-2 far-field modified blocks are replaced
by retained generators. For far-fields, a later ID consumes their contributions directly in
sketch space, as described in \Cref{sec:streaming-id}. Parent construction
regenerates required child interactions and reuses
$D_EG_{E,C}^{\top}$ across sibling row blocks when possible. This removes storage for persistent distance-2 far-field products. For near-fields, the receiving process forms updates from the generators and accumulates them directly into the persistent near-field blocks.

%% file: streaming-id.tex
\section{Sketching-Accelerated Interpolative Decomposition}
\label{sec:streaming-id}
The compression step for the $i$th box $\B$ on level $\ell$ in \Cref{a:fact} requires the ID of the tall submatrix of the currently active matrix $A^\ell_i$, simply denoted by $A_{\F,\B}$. The compression target is $W_\B=\begin{bmatrix}A_{\M(\B),\B}\\K_{\mathrm{proxy},\B}\end{bmatrix} \in\mathbb C^{m_\B\times b_\B}$, where $K_{\mathrm{proxy},\B}$ represents the contribution from the proxy surface \cite{minden2017recursive} whose row count can be negligible, and $A_{\M(\B),\B}$ may have been modified by the Schur complement updates of previously processed boxes. As the number of distance-2 far-field boxes can be as large as $98$ in 3D, compared with $16$ in 2D, the direct ID of $W_\B$ via RRQR can be expensive in both computation and memory ($\bigO(m_\B b_\B^2)$ arithmetic and $\bigO(m_\B b_\B)$ storage). We accelerate the ID step with two optimizations: materialize-then-sketch and stream-and-sketch. The former reduces the size of the target matrix passed to the ID, while the latter avoids storing the full target matrix.

\subsection{Materialize-then-sketch}
\label{subsec:sketch}
\label{subsec:supp-sketch}
Our first optimization retains target assembly but replaces RRQR on $W_\B$
with that on a smaller random projection~\cite{low-rankLiberty} to compute $T$ and $S$ of \eqref{eq:id-approx}:
\begin{equation}
\label{eq:supp-sketch}
 Y_\B=\Omega_\B W_\B\in\mathbb C^{s_\B\times b_\B},
 \qquad \Omega_\B\in\mathbb R^{s_\B\times m_\B}.
\end{equation}
We choose $s_\B=\bigO(b_\B)$, so the reduction is useful when
$m_\B\gg b_\B$. We use a sparse CountSketch variant~\cite{clarkson2017low} with exactly
$\zeta$ nonzeros per column of $\Omega_\B$. Its rows are partitioned into
$\zeta$ nearly equal groups, and each column has one randomly selected
nonzero in each group, with value $\pm1/\sqrt{\zeta}$ chosen by a random
sign; all other entries are zero. In all the experiemnts, we use $s_\B=b_\B$ and $\zeta=4$.
%\chao{Please specify the actual sketch-size rule $s_\B$ (including its multiplier, oversampling, and any cap) and the value of $\zeta$ used in the experiments. Please note any differences across kernels or experimental configurations and how random seeds are selected.}

We first materialize $W_\B$, form $Y_\B$ by sparse projection, and then
compute its ID via RRQR. Formation of $Y_\B$ costs $\bigO(\zeta m_\B b_\B)$, and
the RRQR cost becomes $\bigO(s_\B b_\B k_\B)$.
This reduces QR work while retaining the baseline storage requirement for $W_\B$.

\subsection{Stream-and-sketch}
\label{subsec:stream-target}\label{subsec:sketch-generators}
Our second optimization changes how the same sketch is assembled when generators of \Cref{sec:lazy-schur} are communicated (\texttt{lazy\_schur=1}). Let $W_\B^{(0)}$ contain the original interaction blocks and proxy rows, and let $\Delta W_\B^{(E)}$ contain the updates contributed by a previously eliminated source box $E$, with zero rows where that source contributes no update. By linearity,
\begin{equation}
\label{eq:streamed-sketch}
 Y_\B=\Omega_\B W_\B^{(0)}+\sum_E\Delta Y_\B^{(E)},
 \qquad \Delta Y_\B^{(E)}=\Omega_\B\Delta W_\B^{(E)}.
\end{equation}
The part of original interaction blocks (one block at a time) and proxy rows are sketched directly; the part of the Schur complement updates $\Delta W_\B^{(E)}$ are sketched from their generators without materializing $\Delta W_\B^{(E)}$. 
%\chao{Please confirm that, for modified far-field blocks, the implementation also sketches the original kernel interactions and adds the generator-based update sketches using the same CountSketch indices and signs. Please approve or correct this description and the decomposition in \eqref{eq:streamed-sketch}.} For modified far-field blocks, the update contributions are sketched from their generators without materializing the updated dense blocks. 
From \eqref{eq:generated-pair}, the contribution of source $E$ is
\begin{equation}
  \Delta Y_\B^{(E)}=\sum_A\Omega_A\Delta^{(E)}_{A,\B}=-\!\left(\sum_A\Omega_A G_{E,A}\right)\!\left(D_EG_{E,\B}^{\top}\right).\label{eq:sketch-source-update}
\end{equation}
Here, $\Omega_A$ denotes the columns of $\Omega_\B$ associated with block row $A$, and the sum runs over block rows updated by source $E$.
Therefore, we first sketch the generators $G_{E,A}$ for all row blocks $A$ with $\bigO(\zeta m_{\B} r_E)$ arithmetic, assuming the source box $E$ has redundant size $r_E$. The product of the sketched generators with $D_EG_{E,\B}^{\top}$ is then formed by GEMM with $\bigO(r_E b_\B^2)$ arithmetic. For multiple sources, each generator is sketched and its update contribution accumulated into $Y_\B$ as it arrives. For each block row $A$, the CountSketch indices
and signs are recorded and reused for the original interaction and every source box.

%% file: shared-memory.tex
\section{Shared-Memory Parallelism}
\label{sec:shared-memory}

We now discuss the shared-memory parallelism during local computation, including Lines \ref{l:s_c_1} and \ref{l:s_c_2} of Algorithm \ref{a:color}. In 3D, each refinement multiplies the number of boxes by eight, and coarse levels have few boxes with large ranks. Shared-memory parallelism lets one process per node use all of its cores at fine levels, and splits the work within a box at coarse levels (\Cref{subsec:blas3}). The implementation avoids data races and provides a deterministic accumulation order for Schur complement updates.

\subsection{Color Batches and Deferred Updates}
\label{subsec:thread-parallel}
A \textit{color batch} contains boxes with the same coordinate-parity color $w(\B)$, processed by separate threads. Each distributed boundary color set is one batch; interior boxes form $2^d$ batches. Same-batch boxes are nonadjacent, so one does not change the blocks read by another, but their Schur updates can share neighbor-pair targets. We avoid write races through two phases:
\paragraph{Immediate phase (during the batch)} Each box performs its ID compression and computes the intermediate matrices $\text{temp}_1 = -X_{\S, \R}\, X_{\R, \R}^{-1}$ and $\text{temp}_2 = -X_{\N, \R}\, X_{\R, \R}^{-1}$, which encode the Schur complement contributions to the skeleton and neighbor blocks, respectively. The box only applies the local updates to $A_{\S, \S}$ and $A_{\N, \S}$
from~\eqref{eq:schur}, which do not lead to a race condition. Each box owns its diagonal and persistent near-field blocks, access to explicit or generated distance-2 interactions, and its transformation and pivot factors. It also retains $\text{temp}_2$ until its consumers finish. Neighbor-to-neighbor updates $A_{\N, \N}$ are deferred.

\paragraph{Deferred phase (after the batch)} Once all boxes in the batch have completed the immediate phase, a second pass processes what we call the \emph{affected neighbors} of the batch:
the unordered box pairs $\{\B_k,\B_s\}$ that require explicit Schur complement contributions
$\text{temp}_2$ from at least one box of the batch. Each pair iterates over the batch sources $E\in \N(\B_k)\cap \N(\B_s)$ that contribute to it and computes the Schur complement update $\tilde{A}_{\B_k, \B_s}$ of \eqref{eq:schur} {via GEMM}.
This iteration is executed by one thread and follows a fixed canonical sequence (e.g., sorted by Morton index). The threads are used to parallelize different box pairs $(\B_k,\B_s)$ instead. 

Assigning each affected neighbor pair to one thread removes concurrent writes without locks. That thread pulls contributions from the batch sources in global Morton order, giving a deterministic sequence of GEMMs independent of the thread count and scheduling.

% For lazy interactions, the same source order is used when a later consumer requests the contribution; unneeded dense products are not formed by this pass. This design has three advantages:
% \begin{itemize}
%     \item \textit{Lock-free.} Since each affected neighbor pair $(\B_k,\B_s)$ is processed by exactly one thread and controls its own data, no locks are needed.
%     \item \textit{Deterministic.} The canonical iteration order over sources ensures that floating-point accumulation is reproducible across runs, for a fixed sequence of BLAS calls; the fixed-chunk policy below preserves call shapes as the worker count changes.
%     \item \textit{Efficient.} The deferred GEMM is parallelized over  neighbor pairs $(\B_k,\B_s)$ (which are non-overlapping targets), achieving good thread utilization without contention.
% \end{itemize}

\subsection{BLAS-3 Kernels and Parallelism within a Box}
\label{subsec:blas3}
The generator and sketch representations create dense operations with many
right-hand sides. We implement these operations in their native layout to
reduce temporary storage and expose Level-3 BLAS work. For example, forming
$-X_{\N,\R}D_E^{-1}$ by transposing a tall matrix, applying a left-side solve,
and transposing back requires two full matrix passes and a tall temporary matrix.
Instead, we solve directly from the right. Specifically, for a symmetric-indefinite diagonal block $D_E$, write
$\Pi_E^\top D_E\Pi_E=C_E H_E C_E^\top$, where $\Pi_E$ is a
permutation and $H_E$ contains $1\times1$ and $2\times2$ symmetric pivots.
Then
\[
 X_{\N,\R} D_E^{-1}
 =\biggl(\Bigl(\bigl((X_{\N,\R}\Pi_E)C_E^{-\top}\bigr)H_E^{-1}\Bigr)C_E^{-1}\biggr)\Pi_E^\top.
\]
Bunch--Kaufman factors from LAPACK support the same right-side approach using triangular solves, small pivot solves, and permutations.

The lower levels of the hierarchy usually expose enough boxes to occupy all threads. On coarse
levels, a color batch $\mathcal W$ may contain fewer boxes than the $T_\ell$
workers on a process. A compute-heavy box can then use up to
$\max(1,\lfloor T_\ell/|\mathcal W|\rfloor)$ cooperating
tasks. We can split right-side solves, large GEMMs, and kernel evaluation over the column dimension into fixed-size chunks, so that the BLAS calls are independent of the thread count. The chunk size is chosen to balance BLAS efficiency and task granularity.

% The numerical chunks remain fixed as task concurrency changes: right-side
% solves use 512-row chunks and split GEMMs use 256 output columns per call.
% Task assignment changes which worker executes a chunk, not its BLAS call
% shape or reduction order. We use one BLAS thread per call to avoid nested
% oversubscription. Deferred updates accumulate directly into their target
% with GEMM's $\beta=1$, avoiding a product-sized temporary and a subsequent
% addition pass. Together with canonical source order, these choices control
% rounding variation within a fixed implementation; they do not imply bitwise
% identity between eager and lazy formulations. The current ablation retains
% these kernel optimizations across configurations and does not measure their
% individual speedups.

%% file: analysis.tex
\section{Communication and Memory Analysis}
\label{sec:alg-analysis}
\label{sec:comm-analysis}
\label{s:analysis}
\label{subsec:combined}
\label{subsec:cost-comparison}

We model the communication cost and memory usage of one level of the spatial-coloring algorithm, with and without the on-demand Schur update optimization in \Cref{sec:lazy-schur}. Here, $n_\ell$ denotes the number of boxes and $P$
the number of active processes. A process owns $v^d=n_\ell/P$ boxes and
communicates with at most $q\le 3^d-1$ neighboring processes. The number of boundary boxes on process $p$ is 
\begin{equation}
\label{e:s_ratio}
 |\Gamma_p|=\bigO\!\left((n_\ell/P)^{(d-1)/d}\right).
\end{equation}
For simplicity, we assume that all the Schur updates due to the elimination of a boundary box are needed by a remote process, and the number of near-field boxes and distance-2 far-field boxes per box is $|\N|=3^d-1$ and $|\M|=5^d-3^d$. We further assume that the active block size is $b_\ell$, with skeleton size $k$ (independent of $\ell$ and box ID) and redundant size $r_\ell=b_\ell-k$. At a leaf $b_\ell=n_{leaf}$; for a coarse box assembled from $2^d$ children, $b_\ell\approx2^dk$.

\subsection{Communication cost}
We use the $\alpha$--$\beta$ model for communication analysis. The communication time per process is
\begin{equation}
\label{eq:T-A-block}
 T_{\mathrm{comm},\ell,p}^{(j)}
 \approx 2^d q\alpha+\beta |\Gamma_p|V_{\ell,p}^{(j)}, \qquad j\in\{0,1\}.
\end{equation}
Here, the first term represents the latency cost of $2^d$ rounds of boundary exchanges with $q$ neighboring processes, and $V_{\ell,p}^{(j)}$ is the communication volume per boundary box subject to the Schur update evaluation mode $j\in\{0,1\}$, as described in \Cref{subsec:lazy-modes}. We estimate these volumes as
\begin{align}
 V_{\ell,p}^{(0)}&\simeq \nu b_\ell^2(|\N|(|\N|-1)/2+|\N|),\\
 V_{\ell,p}^{(1)}&\simeq \nu r_\ell^2+\nu r_\ell b_\ell|\N|.\label{eq:mode-volumes}
\end{align}
Here, $\nu$ represents the number of bytes per scalar value. In mode 0, the terms count $|\N|(|\N|-1)/2$ off-diagonal neighbor pairs and $|\N|$ diagonal neighbor blocks, respectively. For an eliminated source box $E$, the two terms in the mode-1 expression of \eqref{eq:mode-volumes} represent the redundant diagonal block $D_E$ and the generator for updating $X_{\N,\N}$, respectively. The reduction of the communication volume from mode 0 to mode 1 is significant, especially when $r_\ell\ll b_\ell$. For example, for a coarse 3D box with $|\N|=26$, $b_\ell\approx8k$, and hence $r_\ell\approx7k$, the ratio is $V_{\ell,p}^{(0)}/V_{\ell,p}^{(1)}\simeq 15$. This payload estimate assumes all neighbor updates are needed remotely; it does not predict a $15\times$ reduction in total communication time, which also includes the latency cost.

\subsection{Stored interaction memory}
We model the per-level storage of off-diagonal interaction blocks and generators. This excludes diagonal blocks, retained factorization data, temporary workspaces, communication buffers, and runtime overhead; it is not a model of total retained memory or per-process peak memory.
In the eager Schur update evaluation mode, both the near-field and distance-2 far-field blocks are stored as dense matrices. In the on-demand evaluation mode, the distance-2 far-field blocks are replaced by the generator representation, which is much smaller in size. The modeled interaction storage per process for the two modes can be estimated as
\begin{align}
\label{eq:interaction-memory}
M_{\ell,p}^{(0)}\simeq\frac{\nu}{2}\frac{n_\ell}{P}|\N|b_\ell^2+\frac{\nu}{2}\frac{n_\ell}{P}|\M|b_\ell^2,\\
M_{\ell,p}^{(1)}\simeq\frac{\nu}{2}\frac{n_\ell}{P}|\N|b_\ell^2+\frac{\nu n_\ell}{P}|\N|b_\ell r_\ell
\end{align}
With $d=3$, $|\N|=26$, $|\M|=98$, $b_\ell=8k$, and $r_\ell=7k$, the ratio is $M_{\ell,p}^{(0)}/M_{\ell,p}^{(1)}\simeq 1.73$, corresponding to a 42\% reduction in modeled interaction storage. Note that per-process peak memory also includes the temporary allocations and can show a different reduction.

%% file: experiments.tex
\section{Experiments} \label{s:experiments}

In this section, we benchmark the optimized spatial-coloring algorithm on three kernels: Laplace,
Helmholtz, and Mat\'ern. All three discretizations satisfy $A=A^\top$. The Laplace and Mat\'ern matrices are real symmetric; Helmholtz is complex
symmetric and uses non-conjugating transposes in the factorization. %The redundant pivot blocks $D_\B$ are factored using LU factorization with partial pivoting for Laplace and Mat\'ern, and Bunch--Kaufman pivoting for Helmholtz (\Cref{subsec:blas3}).
We use regular grids
with $N = n^d$ points to discretize a unit square ($d=2$) or a unit cube ($d=3$), where grid cell $i$ has center coordinate $x_i$ and side length $h=1/n$. The Laplace and Helmholtz runs use $d=3$ and the Mat\'ern runs use
$d=2$. Specifically, the matrix entries for the Laplace kernel are
\begin{equation}
	\label{eq:laplace}
\begin{aligned}
A_{ij} &= \frac{1}{N}\frac{1}{4\pi\|x_i-x_j\|}, i\ne j,
\\
A_{ii} &=
\int_{[-h/2,h/2]^3}
\frac{1}{4\pi\|x\|}\,dx.
\end{aligned}
\end{equation}
The matrix entries for the Helmholtz kernel are
\begin{equation}
	\label{eq:helmholtz}
\begin{aligned}
A_{ij} &= \frac{1}{N}\frac{\exp(\imath \omega\|x_i-x_j\|)}{4\pi\|x_i-x_j\|}, i\ne j,
\\
A_{ii} &=
\int_{[-h/2,h/2]^3}\frac{\exp(\imath \omega\|x\|)}{4\pi\|x\|}dx- \frac{4}{5\omega^2},
\end{aligned}
\end{equation}
corresponding to the following Lippmann--Schwinger equation with excitation function $f(x)$ and wave number $\omega$: 
\begin{equation}
	\label{eq:vie}
\int_{[0,1]^3} \frac{\exp(\imath \omega\|x-y\|)}{4\pi\|x-y\|}u(y)\,dy
-
\frac{4}{5\omega^2}u(x)
=
f(x).
\end{equation}
We define a wave divisor $d_\omega=2\pi n/\omega$ representing the number of grid cells per wavelength along each direction.
Finally, the matrix entries for the Mat\'ern kernel are 
\begin{equation}
\label{eq:matern}
\begin{aligned}
A_{ij}={}&\left(1+\frac{\sqrt{5}\|x_i-x_j\|}{\ell}
+\frac{5\|x_i-x_j\|^2}{3\ell^2}\right)\\
&\quad\cdot\exp\!\left(-\frac{\sqrt{5}\|x_i-x_j\|}{\ell}\right)
+\delta_{ij}\sigma_n^2.
\end{aligned}
\end{equation}
with diagonal regularization $\sigma_n^2$. The Mat\'ern experiments reported here use length scale $\ell=0.1$ and diagonal regularization $\sigma_n^2=10^{-3}$.

The reported factorization time is decomposed as
$
  t_{\mathrm{fact}}
  =
  t_{\mathrm{comp}}
  + t_{\mathrm{comm}},
$
where $t_{\mathrm{comm}}$ is the accumulated boundary data-exchange communication time. The computation component
$t_{\mathrm{comp}}$ (including process idling time) is estimated by subtracting $t_{\mathrm{comm}}$ from the total runtime. %

All experiments were run on the CPU nodes of Perlmutter.  Each node has two AMD EPYC
7763 processors and 512 GB of memory.
Single-threaded BLAS is used in the spatial-coloring algorithm.  The
executables were linked against Google TCMalloc to reduce allocator
fragmentation in these memory-intensive runs.  We also avoid fast-math compiler
optimizations so that the compiler does not apply aggressive floating-point
reorderings that can affect numerical reproducibility. The software environment is GCC 14.3.0 with compiler flags ``-O3 -fopenmp", Cray MPICH 9.1.0, single-threaded OpenBLAS 0.3.30, and LAPACK 3.12.0. 

The optimized 3D timing runs use one MPI process per node, $n_{leaf}=6^3$, and reduction
threshold $n_{\min}=8$. The weak-scaling series with $N/P=96^3$ uses 64 OpenMP workers per process; the larger Laplace series with $N/P=192^3$ uses 128 workers per process. Strong-scaling runs
vary the OpenMP worker count explicitly. The Mat\'ern examples
use 256 processes, eight OpenMP workers per process, and $n_{\min}=4$.

\begin{table}[t]
	\caption{Matvec error $e_{\mathrm{mv}}$ and relative residual error $e_{\mathrm{res}}$ for the spatial-coloring algorithm, tested on the Laplace and Helmholtz kernels. $d_\omega$ denotes the wave divisor for the Helmholtz kernel. %
	}
	\label{t:exp-selected-timings}
	\centering
	\fontsize{8}{10}\selectfont
	\setlength{\tabcolsep}{2.5pt}
	\begin{tabular}{@{}lllrr@{}}
			\toprule
			Problem & $N$ & $d_w$ & $e_{\mathrm{mv}}$ & $e_{\mathrm{res}}$ \\
			\midrule
			Laplace $\epsilon=10^{-3}$ & $96^3$ & -- & $8.05{\times}10^{-4}$ & $3.30{\times}10^{-3}$ \\
			Laplace $\epsilon=10^{-3}$ & $192^3$ & -- & $7.05{\times}10^{-4}$ & $3.93{\times}10^{-3}$ \\
			Laplace $\epsilon=10^{-3}$ & $384^3$ & -- & $2.91{\times}10^{-4}$ & $4.25{\times}10^{-3}$ \\
			Laplace $\epsilon=10^{-3}$ & $768^3$ & -- & $6.81{\times}10^{-4}$ & $4.53{\times}10^{-3}$ \\
			\midrule
			Helmholtz $\epsilon=10^{-4}$ & $96^3$ & 48 & $1.92{\times}10^{-6}$ & $2.68{\times}10^{-6}$ \\
			Helmholtz $\epsilon=10^{-4}$ & $192^3$ & 96 & $6.73{\times}10^{-7}$ & $8.96{\times}10^{-7}$ \\
			Helmholtz $\epsilon=10^{-4}$ & $384^3$ & 192 & $2.12{\times}10^{-7}$ & $3.13{\times}10^{-7}$ \\
			\midrule
			Helmholtz $\epsilon=10^{-3}$ & $96^3$ & 64 & $1.21{\times}10^{-5}$ & $1.66{\times}10^{-5}$ \\
			Helmholtz $\epsilon=10^{-3}$ & $192^3$ & 128 & $4.38{\times}10^{-6}$ & $5.43{\times}10^{-6}$ \\
			Helmholtz $\epsilon=10^{-3}$ & $384^3$ & 256 & $1.45{\times}10^{-6}$ & $1.92{\times}10^{-6}$ \\
			Helmholtz $\epsilon=10^{-3}$ & $768^3$ & 512 & $4.87{\times}10^{-7}$ & $6.55{\times}10^{-7}$ \\
			\bottomrule
		\end{tabular}%
\end{table}

\subsection{Accuracy and Efficacy}
\label{s:exp-acc}

We first test the accuracy and efficacy of the proposed spatial-coloring algorithm on the three kernels. Let $F$ and $F^{-1}$ be the $\mathcal{H}^2$ approximation of $A$ and its inverse, as shown in \eqref{eq:approx-factor} and \eqref{eq:approx-inv}.
Let $u_t$ be a reference solution vector and $u_H=F^{-1}(Au_t)$ be the approximate solution by the spatial-coloring algorithm with a compression tolerance $\epsilon$. The relative matvec error is defined as $e_{\mathrm{mv}}=\|Au_t-Fu_t\|/\|Au_t\|$, and the relative residual is defined as $e_{\mathrm{res}}=\|Au_t-Au_H\|/\|Au_t\|$. %
The matvec associated with $A$ is evaluated using an FFT at machine precision. The solution error is defined as $e_{\mathrm{sol}}=\|u_t-u_H\|/\|u_t\|$.

The errors $e_{\mathrm{mv}}$ and $e_{\mathrm{res}}$ for 3D Laplace with $\epsilon=10^{-3}$, and 3D Helmholtz with $\epsilon=10^{-4}$ and $\epsilon=10^{-3}$ are listed in \Cref{t:exp-selected-timings}. Note that the wave divisor $d_\omega$ is chosen such that the Helmholtz kernel permits not only low numerical ranks, but also well-conditioned $A$. For Laplace, the reported matvec errors are below $10^{-3}$ and residuals range from $3.30\times10^{-3}$ to $4.53\times10^{-3}$. The Helmholtz residuals are smaller, reaching $6.55\times10^{-7}$ at the largest size with $\epsilon=10^{-3}$. %The largest problem size in this table is $N=768^3=452,984,832$ using 512 Perlmutter nodes. %\chao{Do we have $e_{\mathrm{mv}}$ and/or $e_{\mathrm{res}}$ for the largest Laplace run, $N=1536^3=3,623,878,656$? If available, please add the results and evaluation method to this accuracy discussion/table. Otherwise, please distinguish the largest timed run from the largest run with reported accuracy when describing the billion-unknown result.}     

Next, we test the proposed parallel $\mathcal{H}^2$ algorithm on the 2D Mat\'ern kernel. We choose $N=16384^2=268,435,456$, with kernel parameters length scale $\ell=0.1$ and nugget $\sigma_n^2=10^{-3}$ to make $A$ more ill-conditioned \cite{lim2017numerical}. \Cref{t:exp-matern-timings} reports the matvec error and timing breakdown for the spatial-coloring algorithm on 16 Perlmutter nodes. The matvec errors do not reach the desired compression tolerance because proxy surface techniques do not guarantee accuracy for the Mat\'ern kernel \cite{martinsson2025fast}, and we leave the accuracy improvement to future work. That said, the matvec error decreases as the compression tolerance $\epsilon$ is tightened. Due to the ill-conditioning of the matrix, the solution error $e_{\mathrm{sol}}$ obtained by using $F^{-1}$ as a direct solver is large, but the preconditioned conjugate gradient (PCG) method converges in a small number of iterations and reduces the solution error to $e_{\mathrm{sol}}^{\mathrm{cg}}\leq 1.2\times10^{-8}$ for all three runs.

\begin{table}[t]
	\caption{Technical data for the spatial-coloring algorithm on the Mat\'ern kernel.
		These runs use $N=16384^2=268,435,456$, length scale $\ell=0.1$, nugget $\sigma_n^2=10^{-3}$, no
		proxy points, 256 MPI processes on 16 nodes, and reduction threshold $n_{\min}=4$.
		The tree depth is chosen so that each leaf-level box contains $n_{leaf}=8^2=64$
		points. PCG uses $F^{-1}$ as a preconditioner and a convergence tolerance of $10^{-13}$.
		Here $e_{\mathrm{sol}}$ is the direct solution error,
		$n_{\mathrm{cg}}$ is the PCG iteration count, $e_{\mathrm{sol}}^{\mathrm{cg}}$ is the solution error after PCG,
		and $t_{\mathrm{cg}}$ is the time spent in CG iterations including preconditioner application time. All times are in seconds.
		}
	\label{t:exp-matern-timings}
	\centering
	\fontsize{8}{10}\selectfont
	\setlength{\tabcolsep}{1.5pt}
	\begin{tabular}{@{}rrrrrrrrr@{}}
			\toprule
			$\epsilon$ & $t_{\mathrm{fact}}$ & $t_{\mathrm{comm}}$ & $t_{\mathrm{solve}}$ & $e_{\mathrm{mv}}$ & $e_{\mathrm{sol}}$ & $n_{\mathrm{cg}}$ & $e_{\mathrm{sol}}^{\mathrm{cg}}$ & $t_{\mathrm{cg}}$ \\
			\midrule

			$10^{-9}$ & 11.8 & 4.4 & 0.5 & $1.3{\times}10^{-6}$ & $2.5$ & 10 & $1.2{\times}10^{-8}$ & 11.5 \\

			$10^{-10}$ & 12.8 & 4.6 & 0.5 & $9.5{\times}10^{-8}$ & $1.7{\times}10^{-1}$ & 4 & $3.4{\times}10^{-9}$ & 5.2 \\

			$10^{-11}$ & 14.0 & 5.5 & 0.6 & $3.0{\times}10^{-8}$ & $5.7{\times}10^{-2}$ & 4 & $3.4{\times}10^{-9}$ & 5.4 \\

			\bottomrule
		\end{tabular}
\end{table}

\subsection{Ablation Study of Streaming and On-Demand Updates}
\label{s:exp-optimizations}

\Cref{t:exp-sketch-lazy-schur} measures the effects of optimization techniques presented in
\Cref{sec:lazy-schur,sec:streaming-id}. The first row applies direct ID to the materialized target $W_\B$ and
forms Schur products eagerly. The second row enables randomized sketching
while retaining eager Schur updates. The third combines sketching with
on-demand near-field and far-field updates (\texttt{lazy\_schur=1}).
The 3D tests use $N=384^3$, $P=64$, $64$ OpenMP workers and $\epsilon=10^{-3}$; the Mat\'ern test
uses $N=16384^2$, $P=256$, $8$ OpenMP workers and $\epsilon=10^{-11}$.
%\chao{Please confirm the OpenMP worker counts for every ablation configuration in \Cref{t:exp-sketch-lazy-schur}: are all 3D rows run with 64 workers per process, and all Mat\'ern rows with 256 MPI processes and eight workers per process? Please specify any differences between rows.} 

For Helmholtz, dramatic improvements in factorization time and memory have been observed: the complete configuration reduces factorization time from
848.27 s to 99.70 s ($8.5\times$) and peak-process peak memory from 188.15 GiB to
85.65 GiB (54\%). Laplace improves from 714.37 s to 103.72 s ($6.9\times$)
and from 123.24 GiB to 52.91 GiB (57\%). The corresponding exchange-time speedup
ratios are $12.5\times$ and $9.2\times$. These reductions are consistent with the qualitative predictions of \Cref{sec:alg-analysis}. Note that the measured peak-process peak memory includes all resident process memory, whereas the memory model counts only per-level interaction storage. %These are comparisons with the unsketched eager implementation, not with the prior distributed solver.

For Mat\'ern, the combined configuration gives $1.39\times$ lower
factorization time and about 10\% lower RSS. These are reasonable improvements for 2D problems given the analysis in \Cref{sec:alg-analysis}.

\begin{table*}[t]
    \centering
    \caption{Effect of sketching and lazy Schur updates on factorization time, data-exchange communication time, and peak resident memory. The column groups show Helmholtz, Laplace, and Mat\'ern results, respectively. $t_{\mathrm{fact}}$ is the total factorization time, and $t_{\mathrm{comm}}$ is the boundary data exchange communication time. $M_{\max}$ is the maximum \texttt{MaxRSS} in GiB per process. Here \texttt{lazy\_schur}=0 denotes eager updates and \texttt{lazy\_schur}=1 denotes on-demand near-field and far-field updates. \texttt{use\_sketch} indicates whether randomized sketching is used before ID compression.
	}
    \label{t:exp-sketch-lazy-schur}
    \fontsize{8}{10}\selectfont
    \setlength{\tabcolsep}{3pt}
    \begin{tabular}{@{}rr*{3}{rrr}@{}}
        \toprule
        & & \multicolumn{3}{c}{Helmholtz ($N=384^3$)} & \multicolumn{3}{c}{Laplace ($N=384^3$)} & \multicolumn{3}{c}{Mat\'ern ($N=16384^2$)} \\
        \cmidrule(lr){3-5} \cmidrule(lr){6-8} \cmidrule(lr){9-11}
        \texttt{use\_sketch} & \texttt{lazy\_schur}
        & $t_{\mathrm{fact}}$ (s) & $t_{\mathrm{comm}}$ (s) & $M_{\max}$ (GiB)
        & $t_{\mathrm{fact}}$ (s) & $t_{\mathrm{comm}}$ (s) & $M_{\max}$ (GiB)
        & $t_{\mathrm{fact}}$ (s) & $t_{\mathrm{comm}}$ (s) & $M_{\max}$ (GiB) \\
        \midrule
        0 & 0 & 848.27 & 463.19 & 188.15 & 714.37 & 453.22 & 123.24 & 13.94 & 6.23 & 20.21 \\
        1 & 0 & 343.23 & 238.10 & 185.26 & 212.32 & 136.67 & 114.32 & --- & --- & --- \\
% Paper mode 1 denotes the fully lazy configuration; source logs use h2_lazy_schur2.
        1 & 1 & 99.70 & 36.91 & 85.65 & 103.72 & 49.33 & 52.91 & 10.06 & 3.67 & 18.17 \\
        \bottomrule
    \end{tabular}
\end{table*}

\subsection{Parallel Performance}
\label{s:exp-color-scaling}
Next, we investigate the parallelization efficiency of the proposed spatial-coloring algorithm using the 3D Laplace and Helmholtz kernels with a compression tolerance $\epsilon=10^{-3}$.
\subsubsection{Weak scaling}
\begin{figure}[t]
    \centering
    \includegraphics[width=\linewidth]{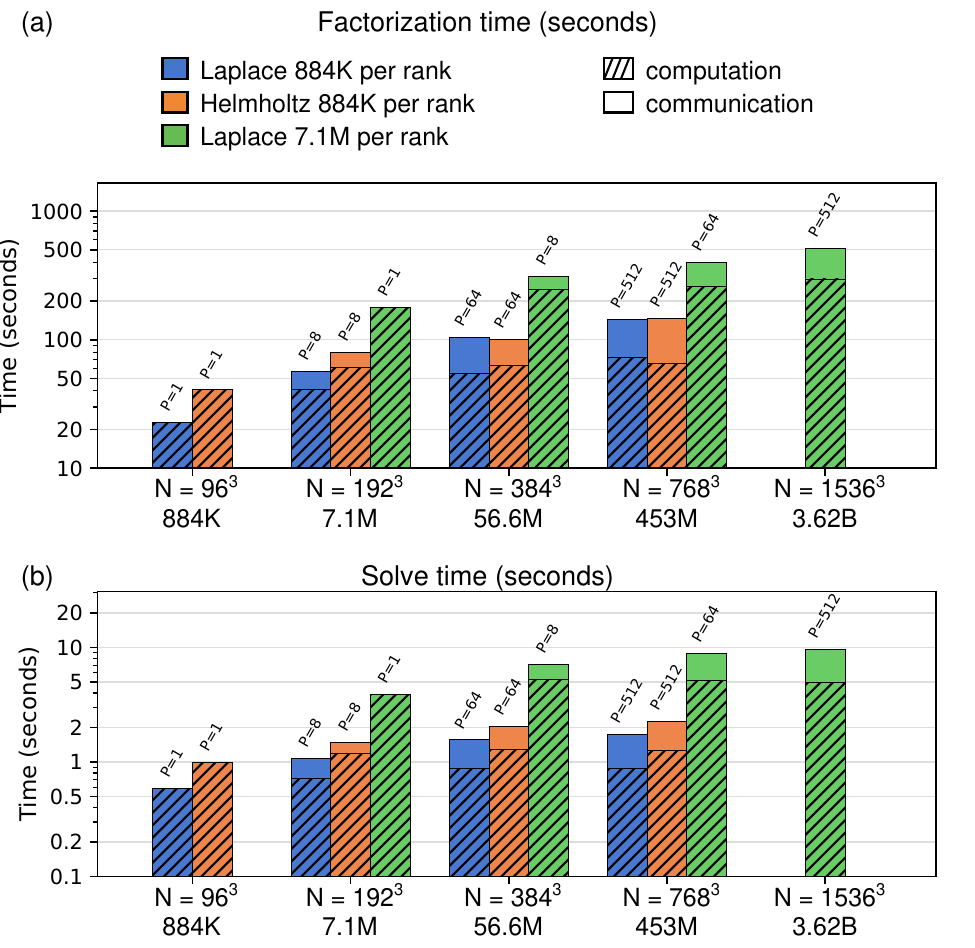}
    \caption{Factorization \textbf{(a)} and solve \textbf{(b)} time breakdowns for the spatial-coloring algorithm on the three-dimensional Laplace and Helmholtz problems, with compression tolerance $\epsilon=10^{-3}$. Each run uses one MPI process per compute node, with $N/P=96^3$ or $192^3$ unknowns per process and up to $P=512$ processes. The panels use separate vertical scales and share the legend. Hatched segments denote computation; solid segments denote communication.}
    \label{f:exp-timing-combined}
    \label{f:exp-laplace-factorization}
    \label{f:exp-helmholtz-e3-factorization}
    \label{f:exp-solve-breakdown}
\end{figure}

\Cref{f:exp-timing-combined} shows the optimized Laplace and Helmholtz runs with 1 MPI process per node. We increase $N$ from $96^3$ to $1536^3$ and $P$ from 1 to 512 while keeping $N/P=96^3=884K$ or $192^3=7.1M$ unknowns per process. To our knowledge, the $N=1536^3=3.62$ billion Laplace run on 512 nodes represents the largest direct solver test case thus far for 3D integral equation kernels. For Laplace with $N/P=96^3$, factorization time grows from 22.79 s to
143.08 s, whereas solve time grows from 0.587 s to 1.751 s. For Laplace with $N/P=192^3$, factorization time grows from 178.1 s to
515 s, whereas solve time grows from 3.9 s to 9.6 s. Similar growth is observed for the Helmholtz runs. The growth is largely due to the increased communication cost (even with the optimization of \Cref{sec:lazy-schur}) and idling time, with only a mild increase in computation time.

For the Laplace runs with 64 OpenMP workers per process and $N/P=96^3=884K$, increasing the problem size from $384^3$ ($P=64$) to $768^3$ ($P=512$) leaves the leaf-level times similar (35.3 s and 31.8 s). The next two levels are also similar, while the main increase occurs in the coarser levels: at $L-\ell=4$, time rises from 14.4 s to 34.5 s, and the larger run adds a level costing 20.5 s.

\Cref{t:exp-level-breakdown-omp128} reports the larger-scale configuration, using 128 OpenMP workers per process and $N/P=192^3=7.1M$, for the $768^3$ ($P=64$) and $1536^3$ ($P=512$) Laplace runs. The increase is largest at $L-\ell=5-6$, with additional communication costs at finer levels. These are separate weak-scaling configurations with different workloads per process and thread counts, rather than a direct comparison of their performance. Both exhibit growth at coarse levels; together with the increased communication times, this supports the explanation that coarse-level workload and reduced concurrency limit weak-scaling efficiency. The coarse-level workload reflects 3D rank growth: at $\epsilon=10^{-3}$, the average skeleton size $k$ rises from $44$ at the leaf level to $97$--$98$ at $L-\ell=3$. At these levels, few boxes remain and each has a large active set. Interestingly, the ranks for 2D integral equation kernel reported in~\cite[Fig.~9]{liang2024on} are smaller and grow more slowly despite the tighter default tolerance $\epsilon=10^{-6}$. %Different tolerances and experimental settings limit this dimensional comparison. 

\begin{table}[t]
\caption{Level-by-level factorization times for the $768^3$ ($P=64$) and $1536^3$ ($P=512$) Laplace runs with 128 OpenMP workers per process. $L-\ell=0$ is the leaf level, $k$ is the average skeleton size, and $\Delta$ is the time difference in seconds (larger minus smaller run). Orange levels performed process reduction. The final coarsening step and root factorization are omitted from the per-level rows but included in the totals.}
\label{t:exp-level-breakdown-omp128}
\centering
\fontsize{8}{10}\selectfont
\setlength{\tabcolsep}{1.5pt}
\begin{tabular}{@{}lr@{\hspace{3pt}}rrr@{\hspace{4pt}}rcr@{\hspace{3pt}}rrr@{\hspace{4pt}}rr@{}}
\toprule
& \multicolumn{5}{c}{$N=768^3$} & & \multicolumn{5}{c}{$N=1536^3$} & \\
\cmidrule(r){2-6}\cmidrule(l){8-12}
$L-\ell$ & $v^d$ & $P$ & $t_{\mathrm{fact}}$ & $t_{\mathrm{comm}}$ & $k$ & & $v^d$ & $P$ & $t_{\mathrm{fact}}$ & $t_{\mathrm{comm}}$ & $k$ & $\Delta$\\
\midrule
0 & 32768 & 64 & 176.2 & 29.3 & 44 & & 32768 & 512 & 195.6 & 37.4 & 44 & $19.3$ \\
1 & 4096 & 64 & 85.3 & 24.7 & 67 & & 4096 & 512 & 96.5 & 30.3 & 67 & $11.2$ \\
2 & 512 & 64 & 45.2 & 13.9 & 86 & & 512 & 512 & 56.9 & 20.7 & 86 & $11.7$ \\
3 & 64 & 64 & 38.9 & 11.2 & 97 & & 64 & 512 & 45.8 & 17.2 & 98 & $6.9$ \\
4 & 8 & 64 & 33.5 & 13.9 & 76 & & 8 & 512 & 36.2 & 14.8 & 76 & $2.6$ \\
\textcolor{commBdy}{\textbf{5}} & 8 & 8 & 21.4 & 7.4 & 81 & & 8 & 64 & 56.6 & 27.1 & 81 & $35.2$ \\
\textcolor{commBdy}{\textbf{6}} & --- & --- & --- & --- & --- & & 8 & 8 & 26.5 & 9.1 & 85 & $26.5$ \\
\midrule
\textbf{Total} & & & \textbf{401.7} & \textbf{100.4} & & & & & \textbf{515.0} & \textbf{156.5} & & $113.3$ \\
\bottomrule
\end{tabular}
\end{table}

\subsubsection{Strong scaling}
In the strong-scaling runs, we fix the Laplace problem size to $N=384^3$ and use 64 MPI processes, with one process per node. For OpenMP thread counts $T=1,4,16,64$, \Cref{t:exp-thread-level-breakdown-both} reports the per-level factorization and communication times. From $T=1$ to $T=64$, the total factorization time decreases from 860 s to 104 s, an approximately $8.3\times$ speedup. The leaf-level time decreases from 461 s to 35 s, while the coarser levels benefit less from additional threads.

\begin{table}[t]
    \caption{Level-by-level breakdown for the strong-scaling runs with $P=64$ for
    the $N=384^3$ Laplace run. Each
    thread count $T$ has two columns: total per-level factorization time ($t_{\mathrm{fact}}$) and communication time ($t_{\mathrm{comm}}$). All times are in seconds. The levels with a reduction of process count, to keep the number of boxes per process $v^d$ sufficiently large, are marked in orange.
    }
    \label{t:exp-thread-level-breakdown-both}
    \centering
    \fontsize{8}{10}\selectfont
       \setlength{\tabcolsep}{1.5pt}
    
        \label{t:exp-color-thread-level-breakdown}
        % Keep the table at an explicit 8-point size; do not scale its contents.
        \begin{tabular}{@{}llr*{4}{rr}@{}}
        \toprule
        $L-\ell$ & $v^d$ & $k$ & \multicolumn{2}{c}{$T=1$} & \multicolumn{2}{c}{$T=4$} & \multicolumn{2}{c}{$T=16$} & \multicolumn{2}{c}{$T=64$} \\
        \cmidrule(lr){4-5} \cmidrule(lr){6-7} \cmidrule(lr){8-9} \cmidrule(lr){10-11}
        & & & $t_{\mathrm{fact}}$ & $t_{\mathrm{comm}}$ & $t_{\mathrm{fact}}$ & $t_{\mathrm{comm}}$ & $t_{\mathrm{fact}}$ & $t_{\mathrm{comm}}$ & $t_{\mathrm{fact}}$ & $t_{\mathrm{comm}}$ \\
        \midrule
        0 & 4096 & 43 & 461 & 41.0 & 122 & 13.7 & 48 & 8.4 & 35 & 7.6 \\
        1 & 512 & 62 & 222 & 42.6 & 58 & 12.0 & 21 & 5.6 & 18 & 5.5 \\
        2 & 64 & 79 & 103 & 37.4 & 29 & 9.4 & 17 & 3.6 & 17 & 4.4 \\
        3 & 8 & 69 & 41 & 23.3 & 20 & 8.8 & 18 & 7.5 & 18 & 7.3 \\
        \textcolor{commBdy}{\textbf{4}} & 8 & 76 & 32 & 9.8 & 17 & 6.0 & 15 & 4.9 & 14 & 4.9 \\
        \midrule
        \textbf{Total} &  &  & \textbf{860} & \textbf{154.1} & \textbf{247} & \textbf{50.0} & \textbf{120} & \textbf{30.0} & \textbf{104} & \textbf{29.7} \\
        \bottomrule
        
        \end{tabular}
\end{table}

%% file: conclusion.tex
\section{Conclusion}

We presented a distributed-memory strong recursive skeletonization algorithm for symmetric $\mathcal{H}^2$ matrices, targeting the larger stencils, box counts, and coarse-level ranks of three-dimensional problems. Spatial box coloring lets every process with eligible boxes participate in each boundary stage, with $2^d$ dependent exchanges per level. On-demand Schur updates and streamed sketching reduce dense storage and communication, while BLAS-3 solves and task splitting improve local execution.

Experiments demonstrate factorization of a 3D Laplace system with 3.62 billion unknowns on 512 nodes and test a 2D Mat\'ern system with 268 million unknowns. To our knowledge, this is the first distributed-memory strong recursive skeletonization factorization demonstrated for three-dimensional volume integral equations at billion-unknown scale. At $N=384^3$, streamed sketching and fully lazy updates reduce factorization time by $6.9\times$ for Laplace and $8.5\times$ for Helmholtz relative to the unsketched, eager configuration, with peak process memory reductions of 57\% and 54\%. Accuracy tests support use as a direct solver or preconditioner.

Future work includes communication-avoiding parallelization, GPU acceleration, and extension to more general kernels and unstructured point distributions.

\section*{Acknowledgments}
This work was supported in part by the U.S. Department of Energy, Office of Science, Office of Advanced Scientific Computing Research, Applied Mathematics Competitive Portfolios program and Scientific Discovery through Advanced Computing (SciDAC) program through the FAST-Math Institute, under contract No. DE-AC02-05CH11231. This research used resources of the National Energy Research Scientific Computing Center (NERSC), a U.S.\ Department of Energy Office of Science User
Facility.